# Stability Analysis of Optimal Adaptive Control using Value Iteration with Approximation Errors

Ali Heydari

*Abstract*—Adaptive optimal control using value iteration initiated from a stabilizing control policy is theoretically analyzed in terms of stability of the system during the learning stage without ignoring the effects of approximation errors. This analysis includes the system operated using any *single/constant* resulting control policy and also using an *evolving/time-varying* control policy. A feature of the presented results is providing estimations of the *region of attraction* so that if the initial condition is within the region, the whole trajectory will remain inside it and hence, the function approximation results remain valid.

## I. INTRODUCTION

Numerous success stories are reported on using adaptive/approximate dynamic programming (ADP) for solving challenging optimal control problems, [1]. Despite the potentials, a shortcoming is that *guarantees* of a suitable performance, needed for using the method in sensitive systems, are yet to be fully established. ADP-based learning algorithms are typically classified as either value iteration (VI) or policy iteration (PI), [1]. PI has the feature that the control under evolution remains stabilizing, [1]. Hence, it is the natural choice for online implementation, i.e., adapting the control 'on the fly.' However, it needs to start with a stabilizing initial control. VI, on the other hand, can be initiated arbitrarily. But, the closed loop system is not guaranteed to be stable during its learning process, if implemented online. It was shown in [2], [3] that in VI also, if the initial control is stabilizing, the control during the learning stage remains stabilizing.

Considering optimal control of discrete-time nonlinear problems with continuous state and action spaces and undiscounted cost functions using VI, which is the subject of this work, learning convergence was investigated by different researchers including [4]–[8]. All these convergence analyses are based on the assumption of perfect function reconstruction, i.e., no error in the function approximation. However, the approximation errors exist almost in every application when the system is nonlinear or when cost function terms are non-quadratic and nonlinear. What makes their presence potentially problematic is the fact that the errors *propagate* throughout the iterations. Hence, regardless of how small they are, their accumulated effect can become significant, leading to a phenomenon similar to *resonance*, and hence, unreliability of the result, [9].

Analyzing VI under the presence of approximation errors, i.e. *approximate VI* (AVI), is an open research problem with a few published results, including [9]–[13]. Refs. [10]–[12] investigated problems with *discounted* cost functions and the results are solely valid for such problems, as the 'forgetting'

A. Heydari is an Assistant Professor of Mechanical Engineering with the Southern Methodist University, Dallas, TX, email: aheydari@smu.edu.

nature of discounted problems is the backbone of the development of the error bounds. On the other hand, the results in [13] provided error analyses but with assumptions whose verification is not straight forward. For example, the approximation error between the exact and approximate functions must be possible to be written in a *multiplicative* form in [13], instead of an *additive* form, as done in this study. As for [9], it investigated AVI with *arbitrary* initial guesses and *after* the training stage, i.e., it extended the result in [7] to AVI. The current study, however, investigates AVI initiated using a *stabilizing* guess and *during* the training phase. Based on this background, theoretical analysis of the consequences of the approximation errors on the results is of great interest to the ADP researchers and practitioners. Finally, Ref. [14] may be mentioned for an analysis of PI with approximation errors.

In terms of contributions of this study, initially, stability of the system operated using any single control policy generated using AVI is investigated. Afterwards, the legitimate concern that any ADP result is valid only when the state trajectory remains within the region for which the controller is trained, is addressed. This concern is resolved through establishing an *estimation of the region of attraction* (EROA) [15] for the controller, such that as long as the initial condition of the system is within the region, the entire trajectory is guaranteed to remain in the region. Then, these results are extended to the case of applying an evolving control policy, i.e., updating the policy which is being used for control, on the fly.

The rest of this study is organized as follows. The optimal control problem is formulated in Section II and the value iteration based solution is revisited in Section III. Section IV presents the main results, i.e., theoretical analyses of the learning scheme. Afterward, some numerical results are given in Section V, followed by concluding remarks in Section VI.

## II. PROBLEM FORMULATION

Let the system subject to control be given by discrete-time nonlinear dynamics

$$x_{k+1} = f(x_k, u_k), k \in \mathbb{N}, \quad (1)$$

where $f: \mathbb{R}^n \times \mathbb{R}^m \to \mathbb{R}^n$ is a continuous function versus its both inputs, i.e., the state and control vectors, $x$ and $u$, respectively, with $f(0,0) = 0$. The sets of real numbers and non-negative integers are denoted with $\mathbb{R}$ and $\mathbb{N}$, respectively. Positive integers $n$ and $m$ denote the dimensions of the continuous state and control spaces. Finally, sub-index $k$ represents the discrete time index. The cost function is given by

$$J = \sum_{k=0}^{\infty} U(x_k, u_k), \quad (2)$$

where *utility function* $U(\cdot, \cdot)$ is of form $U(x_k, u_k) := Q(x_k) + R(u_k)$ for continuous and positive definite functions $Q : \mathbb{R}^n \to \mathbb{R}_+$ and $R : \mathbb{R}^m \to \mathbb{R}_+$. Set $\mathbb{R}_+$ denotes the non-negative reals. Let *control policy* $h : \mathbb{R}^n \to \mathbb{R}^m$ be used for feedback control calculation, i.e., $u_k = h(x_k)$. The objective is finding the *optimal control policy*, denoted with $h^*(\cdot)$, that is the policy using which cost function (2) is minimized, subject to (1). In online learning, this process is done through selecting an initial control policy and *updating/adapting* it, until it converges to the optimal control policy.

**Definition 1.** *A control policy $h(\cdot)$ is defined to be admissible within a connected and compact set $\Omega \subset \mathbb{R}^n$ containing the origin, if a) it is a continuous function of $x$ in the set with $h(0) = 0$ and b) its 'cost-to-go' or 'value function', denoted with $V_h : \mathbb{R}^n \to \mathbb{R}_+$ and defined by*

$$V_h(x_0) = \sum_{k=0}^{\infty} U\big(x_k^h, h(x_k^h)\big), \quad (3)$$

*is continuous in $\Omega$. In Eq. (3) one has $x_k^h := f\big(x_{k-1}^h, h(x_{k-1}^h)\big), \forall k \in \mathbb{N} - \{0\}$, and $x_0^h := x_0$. In other words, $x_k^h$ denotes the kth element on the state history initiated from $x_0$ and propagated using control policy $h(\cdot)$.*

The main difference between the defined admissibility and the ones typically utilized in the ADP literature, including [5], is the assumption of continuity of the value function. The continuity is required for *uniform* approximation of the function using parametric function approximators, [16]. Also, continuous functions are bounded in compact sets [17], hence, the continuity of the value function leads to its boundedness. This is an essential requirement for an admissible control.

The following assumption applies to the entire results presented in this study and guarantees that there is no state vector in $\Omega$ for which the value function associated with the *optimal* control policy is unbounded.

**Assumption 1.** *There exists at least one admissible control policy for the given system within $\Omega$.*

### III. VALUE ITERATION INITIATED USING STABILIZING CONTROL

Defining the *optimal value function*, as the value function associated with the optimal control policy and denoting it with $V^*(\cdot)$, Bellman equation [18], given below, provides the solution to the problem

$$h^*(x) \in \arg\min_{u \in \mathbb{R}^m} \Big(U(x,u) + V^*\big(f(x,u)\big)\Big), \quad (4)$$

$$V^*(x) = \min_{u \in \mathbb{R}^m} \Big(U(x,u) + V^*\big(f(x,u)\big)\Big). \quad (5)$$

It is worth mentioning that the minimizing $u$ in (4) may not be unique and notation $\in$ used here (motivated by [8]) allows selecting any of the minimizers. Due to the *curse of dimensionality* [18], however, the proposed solution is computationally impracticable for general nonlinear systems. ADP utilizes the idea of *approximating* the optimal value function, using either look-up tables or function approximators, e.g., artificial neural networks, for remedying the problem. The value function approximator is typically called the *critic* in the ADP literature. The approximation is performed over a *compact* and *connected* set containing the origin, called the *region of interest*. This region, denoted with $\Omega$, has to be selected based on the specific problem at hand and the expected range of states to be visited during operation of the system. It should be noted that the ADP-based results are valid only if the entire state trajectory initiated from the initial state vector remains within $\Omega$.

Approximation of optimal value function can be done using VI. The process starts with an initial guess $V^0(\cdot)$ and iterates through the *policy update equation* given by

$$h^i(x) \in \arg\min_{u \in \mathbb{R}^m} \Big(U(x,u) + V^i\big(f(x,u)\big)\Big), \forall x \in \Omega, \quad (6)$$

and the *value update equation*

$$V^{i+1}(x) = U\big(x, h^i(x)\big) + V^i\Big(f\big(x, h^i(x)\big)\Big), \forall x \in \Omega, \quad (7)$$

or equivalently

$$V^{i+1}(x) = \min_{u \in \mathbb{R}^m} \Big(U(x,u) + V^i\big(f(x,u)\big)\Big), \forall x \in \Omega, \quad (8)$$

for $i = 0, 1, ...$ until the iterations converge. If the iterations converge to the optimal value function, i.e., if $V^i(\cdot) \to V^*(\cdot)$ as $i \to \infty$, the resulting $V^*(\cdot)$ can be used in (4) for finding the optimal policy.

A critical concern, when using VI for online learning, is the stability of the system under *immature* control policies, i.e., $h^i(\cdot)$s before the convergence of the solution. This stability is not guaranteed. However, it was shown that if the initial guess $V^0(\cdot)$ is selected as the value function of an admissible policy, then, the immature policies remain stabilizing, [2]. Denoting the initial admissible policy with $h^{-1}(.)$ (for notational compatibility,) its value function, denoted with $V^0(.)$, can be calculated through (3), or equivalently, by solving

$$V^0(x) = U\big(x, h^{-1}(x)\big) + V^0\Big(f\big(x, h^{-1}(x)\big)\Big), \forall x \in \Omega. \quad (9)$$

### IV. STABILITY ANALYSIS UNDER APPROXIMATE VALUE ITERATION

*Exact* reconstruction of the right hand side of Eq. (8) is generally impossible except for simple problems. In practice, one utilizes parametric function approximators for this purpose. Such approximators introduce *approximation errors* into the process. The approximation errors replace (8) with

$$\hat{V}^{i+1}(x) = \min_{u \in \mathbb{R}^m} \Big(U(x,u) + \hat{V}^i\big(f(x,u)\big)\Big) + \epsilon^i(x), \forall x \in \Omega, \quad (10)$$

where the *approximate value function* at the $i$th iteration is denoted with $\hat{V}^i(\cdot)$ and the approximation error 'at this iteration' is denoted with $\epsilon^i(\cdot)$. Note that the value function in the right hand side of Eq. (10) is also an approximate quantity, generated from the previous iteration. Therefore, $\epsilon^i(\cdot)$ 'is not' the difference between the exact value function $V^{i+1}(\cdot)$ and the approximate one $\hat{V}^{i+1}(\cdot)$. The difference between $V^{i+1}(\cdot)$ and $\hat{V}^{i+1}(\cdot)$ includes the *accumulation* of the $\epsilon^i(\cdot)$s throughout the conducted $i$ iterations. When $\epsilon^i(\cdot) \neq 0$, convergence of the approximate VI (AVI) does not follow from Ref. [2].

Before proceeding to the analysis with approximation errors, it is worth mentioning that one typically trains a control approximator (*actor*) to approximate the solution to the minimization problem given by (4) based on the value function resulting from VI. The actor, will hence, lead to *another* approximation error term in the process, regardless of whether the value function reconstruction is exact or approximate. However, the effect of the actor's approximation error can be removed from both the convergence analysis of AVI and the stability analysis of the system during AVI. Because, the control will be directly calculated from the minimization of the right hand side of Eq. (10) in *online* and *adaptive* optimal control and applied on the system. In other words, even though the actor will be updated simultaneously along with the critic in online learning, the critic training and the operation of the system are independent of approximation accuracy of the actor. Once the learning is concluded (and if it is concluded,) the operation of the system could be based on the control resulting from the trained actor, hence, the approximation error of actor can affect the stability of the system at that stage. The stability analysis *after* conclusion of AVI is investigated by the author in [9]. The focus in this work is analyzing AVI *during* the learning process. Finally, it is worth mentioning that the AVI analyses presented next are valid for the case of *offline* training as well. The reason is, in offline training also, the actor training can be postponed till after the convergence of the value iterations, as detailed in [7]. Therefore, the actor approximation errors will not affect the convergence of AVI. Considering the abovementioned comment and denoting the minimizer of the right hand side of Eq. (10) by $\hat{h}^i(\cdot)$, one has

$$\hat{h}^i(x) \in \arg\min_{u \in \mathbb{R}^m} \Big(U(x,u) + \hat{V}^i\big(f(x,u)\big)\Big), \forall x \in \Omega, \quad (11)$$

therefore, Eq. (10) can be written as

$$\hat{V}^{i+1}(x) = U\big(x, \hat{h}^i(x)\big) + \hat{V}^i\big(f(x, \hat{h}^i(x))\big) + \epsilon^i(x), \forall x \in \Omega. \quad (12)$$

Given initiation of exact VI (i.e., the VI without approximation errors) using the *exact* value function of an admissible control policy, let the AVI be initiated using an *approximation* of the value function of an admissible control policy. Let the approximation error be denoted with $\epsilon^{-1}(\cdot)$. In other words instead of the exact value function $V^0(\cdot)$ given by (9), one initiates the iterations using the approximate value function $\hat{V}^0(\cdot)$ which satisfies

$$\hat{V}^0(x) = U(x, h^{-1}(x)) + \hat{V}^0\big(f(x, h^{-1}(x))\big) \\ + \epsilon^{-1}(x), \forall x \in \Omega. \quad (13)$$

It should be noted that the difference between $V^0(x)$ and $\hat{V}^0(x)$ is not simply given by $\epsilon^{-1}(x)$. This would have been the case 'if' instead of $\hat{V}^0(\cdot)$, function $V^0(\cdot)$ was used in the right hand side of (13). As a matter of fact, Eq. (13) leads to

$$\hat{V}^0(x_0) = \sum_{k=0}^{\infty} \Big(U\big(x_k^{h^{-1}}, h^{-1}(x_k^{h^{-1}})\big) + \epsilon^{-1}(x_k^{h^{-1}})\Big), \quad (14)$$

while

$$V^0(x_0) = \sum_{k=0}^{\infty} U\big(x_k^{h^{-1}}, h^{-1}(x_k^{h^{-1}})\big). \quad (15)$$

Hence, if $|\epsilon^{-1}(x)| \le cU(x,0), \forall x$, for some $c \in [0,1)$ which leads to

$$(1-c)U(x,u) \le U(x,u) + \epsilon^{-1}(x) \le (1+c)U(x,u), \forall x, \forall u, \quad (16)$$

the relation between $V^0(x)$ and $\hat{V}^0(x)$ can be established as

$$(1-c)V^0(x) \le \hat{V}^0(x) \le (1+c)V^0(x), \forall x. \quad (17)$$

Assuming an upper bound for the approximation error $\epsilon^i(x)$ the results given by Theorem 1 (adapted from [9]) can be obtained. They are in terms of boundedness of sequence $\{\hat{V}^i(x)\}_{i=0}^{\infty}$ resulting from the AVI and its relation versus the optimal value function. This boundedness will later be used for stability analysis. Before that, an assumption is made.

**Assumption 2.** *The approximate value iteration given by (12) is initiated using the approximate value function of an admissible policy (denoted with $h^{-1}(\cdot)$ in this study) and conducted using a continuous function approximator with the bounded approximation error $|\epsilon^i(x)| \le cU(x,0), \forall i \in \mathbb{N} \cup \{-1\}, \forall x \in \Omega$, for some $c \in [0,1)$.*

**Theorem 1.** *Let $\{\overline{V}^i(x)\}_{i=0}^{\infty}$ and $\{\underline{V}^i(x)\}_{i=0}^{\infty}$ be defined as sequences of functions generated using the exact value iterations corresponding to cost functions $\underline{J}$ and $\overline{J}$, respectively, where*

$$\underline{J} = \sum_{k=0}^{\infty} \Big(U(x_k, u_k) - cU(x_k, 0)\Big), \quad (18)$$

$$\overline{J} = \sum_{k=0}^{\infty} \Big(U(x_k, u_k) + cU(x_k, 0)\Big). \quad (19)$$

*Also, let the iterations be initiated using the value functions (based on the respective cost functions) of an admissible policy. If Assumption 2 holds and the approximate value iteration is initiated using some $\hat{V}^0(x)$ which satisfies $\underline{V}^0(x) \le \hat{V}^0(x) \le \overline{V}^0(x), \forall x \in \Omega$, then $\underline{V}^i(x) \le \hat{V}^i(x) \le \overline{V}^i(x), \forall i \in \mathbb{N}, \forall x \in \Omega$.*

*Proof*: The proof (which resembles Relaxed Dynamic Programming presented in [4]) provided in [9] applies to the case of initiating AVI using an admissible policy as well. □

Sequences $\{\overline{V}^i(x)\}_{i=0}^{\infty}$ and $\{\underline{V}^i(x)\}_{i=0}^{\infty}$ converge as shown in [2], when initiated using value functions (defined based on the respective cost function) of some admissible controls. Therefore, considering Theorem 1, the boundedness of the AVI-based results follows, when Assumption 2 holds.

The actual convergence as well as stability of the system operated under AVI are much more challenging compared to the respective analyses in exact VI. The reason is, the presence of approximation errors cancels the monotonicity feature of value functions during the learning stage, [2]. As long as the boundedness of the functions during AVI is guaranteed in a neighborhood of the optimal value function (Theorem 1) where the neighborhood shrinks if the approximation error decreases, the actual convergence of the iteration may not be of critical importance in implementing AVI. But, stability of the system operated under AVI-based results is critical.

The following lemma develops a 'semi-monotonicity' for AVI to be used later for deriving desired stability results.

**Lemma 1.** *Let $\check{V}^i(x_0) := \sum_{k=0}^i U(\hat{x}_k^{*,i}, 0), \forall x_0 \in \Omega, \forall i \in \mathbb{N}$,*

where $\hat{x}_k^{*,i} := f\big(\hat{x}_{k-1}^{*,i}, \hat{h}^{i-k}(\hat{x}_{k-1}^{*,i})\big)$ and $\hat{x}_0^{*,i} := x_0, \forall i, \forall k \in \mathbb{N} - \{0\}$. *If Assumption 2 holds, then,*

$$\hat{V}^{i+1}(x) \leq \hat{V}^i(x) + 2c\breve{V}^i(x), \forall x \in \Omega. \quad (20)$$

*Proof*: Initially note that $\hat{x}_k^{*,i}$, in the statement of the lemma, is the $k$th state vector on the state history initiated from $x_0$ and propagated by applying control policy $\hat{h}^{(i-1)-\bar{k}}(\cdot)$ at time $\bar{k}, 0 \leq \bar{k} \leq i$. The first iteration of AVI leads to

$$\hat{V}^1(x) = \min_{u \in \mathbb{R}^m} \Big(U(x,u) + \hat{V}^0\big(f(x,u)\big)\Big) + \epsilon^0(x). \quad (21)$$

One has $\hat{V}^1(x) - \epsilon^0(x) \leq \hat{V}^0(x) - \epsilon^{-1}(x), \forall x$, because, per (21), $\hat{V}^1(x)$ is resulted from a minimization, as opposed to using a given policy $h^{-1}(\cdot)$ in (13). The foregoing inequality along with $-cU(x,0) \leq -\epsilon^i(x) \leq cU(x,0), \forall i$, lead to

$$\hat{V}^1(x) \leq \hat{V}^0(x) + 2cU(x,0), \forall x, \quad (22)$$

which confirms that inequality (20) holds for $i = 0$. Now, assume that

$$\hat{V}^i(x) \leq \hat{V}^{i-1}(x) + 2c\breve{V}^{i-1}(x), \forall x \in \Omega, \quad (23)$$

if this assumption leads to (20), the proof will be complete by induction. Let

$$\hat{\mathcal{V}}(x) := U(x, \hat{h}^{i-1}(x)) + \hat{V}^i\Big(f\big(x, \hat{h}^{i-1}(x)\big)\Big) + \epsilon^i(x). \quad (24)$$

Since the minimizer of the right hand side of (10) is $\hat{h}^i(\cdot)$, per (11), and not $\hat{h}^{i-1}(\cdot)$, one has

$$\hat{V}^{i+1}(x) \leq \hat{\mathcal{V}}(x), \forall x \in \Omega. \quad (25)$$

On the other hand, comparing $\hat{V}^i(x_0)$, given by replacing $i$ with $i-1$ in (12), with (24) and considering (23) one has

$$\hat{\mathcal{V}}(x) - \epsilon^i(x) \leq \hat{V}^i(x) - \epsilon^{i-1}(x) \\ + 2c\breve{V}^{i-1}\Big(f\big(x, \hat{h}^{i-1}(x)\big)\Big), \forall x \in \Omega, \quad (26)$$

which given $-cU(x,0) \leq -\epsilon^i(x) \leq cU(x,0), \forall i$, leads to

$$\hat{\mathcal{V}}(x) \leq \hat{V}^i(x) + 2cU(x,0) + 2c\breve{V}^{i-1}\Big(f\big(x, \hat{h}^{i-1}(x)\big)\Big), \quad (27)$$

and along with (25) gives

$$\hat{V}^{i+1}(x) \leq \hat{V}^i(x) + 2cU(x,0) \\ + 2c\breve{V}^{i-1}\Big(f\big(x, \hat{h}^{i-1}(x)\big)\Big), \forall x \in \Omega. \quad (28)$$

The next step is showing that

$$U(x_0, 0) + \breve{V}^{i-1}\Big(f\big(x_0, \hat{h}^{i-1}(x_0)\big)\Big) = \breve{V}^i(x_0). \quad (29)$$

Note that, $\breve{V}^i(x_0)$ is the result of evaluating a finite sum of $U(x_k, 0)$'s along a 'trajectory' initiated from $x_0$. So, in order to show that (29) holds, it suffices to show that the summations in both sides of (29), each having $i+1$ summands, are along the same trajectory.

The first summand in $\breve{V}^i(x_0)$ is $U(x_0, 0)$ which is matched by the same term existing in the left hand side of (29). The second summand of $\breve{V}^i(x_0)$ is $U(\cdot, 0)$ evaluated at $\hat{x}_1^{*,i} = f\big(\hat{x}_0^{*,i}, \hat{h}^{i-1}(\hat{x}_0^{*,i})\big)$. The first summand of $\breve{V}^{i-1}\big(f(x_0, \hat{h}^{i-1}(x_0))\big)$ is $U(\cdot, 0)$ evaluated at $x_1 := f\big(x_0, \hat{h}^{i-1}(x_0)\big)$. Since $\hat{x}_0^{*,i} = x_0$, one has $\hat{x}_1^{*,i} = x_1$, hence, the second summand of $\breve{V}^i(x_0)$ also will be matched by a term in the left hand side of (29). Similarly, the third summand of $\breve{V}^i(x_0)$ is $U(\cdot, 0)$ evaluated at $\hat{x}_2^{*,i} = f\big(\hat{x}_1^{*,i}, \hat{h}^{i-2}(\hat{x}_1^{*,i})\big)$.

The second summand of $\breve{V}^{i-1}\big(f(x_0, \hat{h}^{i-1}(x_0))\big)$ is $U(\cdot, 0)$ evaluated at $x_2 := f\big(x_1, \hat{h}^{i-2}(x_1)\big)$, by definition. Since $\hat{x}_1^{*,i} = x_1$ one has $\hat{x}_2^{*,i} = x_2$. Repeating this argument it is seen that the trajectories are identical and hence, (29) holds, which along with (28) proves (20). $\square$

Next, the stability is investigated and an EROA is established. Before that, the term EROA is formally defined.

**Definition 2.** *An estimation of the region of attraction (EROA) for the closed loop system $x_{k+1} = f(x_k, h(x_k))$ is given by $\mathcal{B} \subset \mathbb{R}^n$ if any state trajectory of the system initiated inside $\mathcal{B}$ is defined and converges to the origin as $k \to \infty$, [15].*

**Theorem 2.** *Let Assumption 2 hold. For every given $i \in \mathbb{N}$, control policy $\hat{h}^i(\cdot)$ asymptotically stabilizes the system about the origin if $c$ is such that*

$$0 \leq c < 1 + 2\gamma - \sqrt{4\gamma^2 + 4\gamma}, \quad (30)$$

*for some $\gamma \in \mathbb{R}_+$ which satisfies $V^0(x) \leq \gamma U(x,0), \forall x \in \Omega$ where $V^0(\cdot)$ denotes the exact value function of $h^{-1}(\cdot)$. Moreover, let the compact region $\hat{\mathcal{B}}_r^i$ for any $r \in \mathbb{R}_+$ be defined as $\hat{\mathcal{B}}_r^i := \{x \in \mathbb{R}^n : \hat{V}^i(x) \leq r\}$ and let $\bar{r}^i > 0$ be (possibly the greatest $r$) such that $\hat{\mathcal{B}}_{\bar{r}^i}^i \subset \Omega$. Then, $\hat{\mathcal{B}}_{\bar{r}^i}^i$ will be an estimation of the region of attraction for the system.*

*Proof*: The idea is using $\hat{V}^i(\cdot)$ as a candidate Lyapunov function to prove the claim. The lower and upper boundedness of the function, established in Theorem 1, guarantees the positive definiteness of the function and the continuity of the parametric function approximator guarantees its continuity. The objective is showing negative definiteness of $\Delta \hat{V}^i(x) := \hat{V}^i(f(x, \hat{h}^i(x))) - \hat{V}^i(x)$. By Eq. (12)

$$\hat{V}^i\Big(f\big(x, \hat{h}^i(x)\big)\Big) - \hat{V}^{i+1}(x) = -U\big(x, \hat{h}^i(x)\big) - \epsilon^i(x). \quad (31)$$

Lemma 1 and inequality (20) may be used in the foregoing equation to replace $\hat{V}^{i+1}(x)$ with $\hat{V}^i(x)$ in its left hand side. Before that, let us show that

$$\hat{V}^i(\hat{x}_0^{*,i}) = \\ \sum_{k=0}^{i-1} \Big(U\big(\hat{x}_k^{*,i}, \hat{h}^{(i-1)-k}(\hat{x}_k^{*,i})\big) + \epsilon^{(i-1)-k}(\hat{x}_k^{*,i})\Big) + \quad (32) \\ \hat{V}^0(\hat{x}_i^{*,i}), \forall \hat{x}_0^{*,i} \in \Omega, \forall i \in \mathbb{N} - \{0\},$$

where $\hat{x}_k^{*,i} := f\big(\hat{x}_{k-1}^{*,i}, \hat{h}^{i-k}(\hat{x}_{k-1}^{*,i})\big)$ and $\hat{x}_0^{*,i} := x_0, \forall i, \forall k \in \mathbb{N} - \{0\}$. This can be shown by induction. Eq. (32) for $i = 1$ leads to

$$\hat{V}^1(\hat{x}_0^{*,1}) = U\big(\hat{x}_0^{*,1}, \hat{h}^0(\hat{x}_0^{*,1})\big) + \epsilon^0(\hat{x}_0^{*,1}) + \hat{V}^0(\hat{x}_1^{*,1}), \quad (33)$$

which holds per (12). Now, let us assume Eq. (32) holds for some $i$. The induction is complete if it can be shown that it holds for $i+1$ also, that is

$$\hat{V}^{i+1}(\hat{x}_0^{*,i+1}) = \\ \sum_{k=0}^{i} \Big(U\big(\hat{x}_k^{*,i+1}, \hat{h}^{i-k}(\hat{x}_k^{*,i+1})\big) + \epsilon^{i-k}(\hat{x}_k^{*,i+1})\Big) + \quad (34) \\ \hat{V}^0(\hat{x}_{i+1}^{*,i+1}), \forall \hat{x}_0^{*,i+1} \in \Omega,$$

Given (12),
$$\hat{V}^{i+1}(\hat{x}_0^{*,i+1}) = U(\hat{x}_0^{*,i+1}, \hat{h}^i(\hat{x}_0^{*,i+1})) + \epsilon^i(\hat{x}_0^{*,i+1}) + \hat{V}^i\Big(f\big(\hat{x}_0^{*,i+1}, \hat{h}^i(\hat{x}_0^{*,i+1})\big)\Big), \forall \hat{x}_0^{*,i+1} \in \Omega. \quad (35)$$

Evaluating (32) at $\hat{x}_1^{*,i+1} = f(\hat{x}_0^{*,i+1}, \hat{h}^i(\hat{x}_0^{*,i+1}))$ and using the resulting right hand side for replacing $\hat{V}^i(\cdot)$ in the right hand side of (35), Eq. (34) can be shown to be resulted. The first two summands in (35) are $U(\hat{x}_0^{*,i+1}, \hat{h}^i(\hat{x}_0^{*,i+1})) + \epsilon^i(\hat{x}_0^{*,i+1})$ which are matched by the first set of summands in the summation in (34). The second set of summands in (35) is $U(\hat{x}_1^{*,i+1}, \hat{h}^{i-1}(\hat{x}_1^{*,i+1})) + \epsilon^{i-1}(\hat{x}_1^{*,i+1})$, given evaluating (32) at $\hat{x}_1^{*,i+1}$ when replacing $\hat{V}^i(\cdot)$. This set is also matched by the second set of summands in the summation in (34). The third set of summands in (35) is $U(\cdot, \hat{h}^{i-2}(\cdot)) + \epsilon^{i-2}(\cdot)$ evaluated at $\hat{x}_1^{*,i} = f(\hat{x}_0^{*,i}, \hat{h}^{i-1}(\hat{x}_0^{*,i}))$. But, the third set in (34) is $U(\cdot, \hat{h}^{i-2}(\cdot)) + \epsilon^{i-2}(\cdot)$ evaluated at $\hat{x}_2^{*,i+1} = f(\hat{x}_1^{*,i+1}, \hat{h}^{i-1}(\hat{x}_1^{*,i+1}))$. Given $\hat{x}_0^{*,i} = \hat{x}_1^{*,i+1}$ which follows from evaluating (32) at $\hat{x}_1^{*,i+1}$ and using the resulting right hand side for replacing $\hat{V}^i(\cdot)$ in (35), one has $\hat{x}_1^{*,i} = \hat{x}_2^{*,i+1}$. Therefore, the third sets also match. This process can continue to show that all the summands are the same functions evaluated at the same states. Therefore, (32) holds.

Considering $(1-c)U(x,0) \le U(x,u) + \epsilon^i(x), \forall x, \forall u, \forall i$, and $(1-c)U(\hat{x}_i^{*,i}, 0) \le \hat{V}^0(\hat{x}_i^{*,i})$, comparing (32) with $\breve{V}(x_0)$ defined in Lemma 1, one has

$$(1-c)\breve{V}^i(x) \le \hat{V}^i(x), \forall x \in \Omega. \quad (36)$$

Moreover, by Theorem 1, one has $\hat{V}^i(x) \le \overline{V}^i(x), \forall x$, if $\overline{V}^i(\cdot)$ is generated using the value function of $h^{-1}(\cdot)$ as the initial guess. Moreover, $\overline{V}^i(x) \le \overline{V}^0(x), \forall x$, per monotonicity of exact VI, [2]. Therefore, (36) leads to

$$\breve{V}^i(x) \le \frac{1}{1-c}\overline{V}^0(x), \forall x \in \Omega, \forall i \in \mathbb{N}. \quad (37)$$

The interesting point about inequality (37) is showing the boundedness of the left hand side for any $i$. This boundedness can be used to show the stability of the state trajectory $\hat{x}_k^{*,i}$, given the fact that the left hand side is composed of a partial sum over this trajectory. But, this is neither the trajectory whose stability is under investigation in this theorem, nor the one to be investigated later in Theorem 3. This boundedness, however, along with the results from Lemma 1 lead to

$$\hat{V}^{i+1}(x) \le \hat{V}^i(x) + \frac{2c}{1-c}\overline{V}^0(x), \forall x \in \Omega, \forall i \in \mathbb{N}. \quad (38)$$

Utilizing (38) in (31) leads to

$$\hat{V}^i\big(f(x, \hat{h}^i(x))\big) - \hat{V}^i(x) \le -U(x, \hat{h}^i(x)) - \epsilon^i(x) + \frac{2c}{1-c}\overline{V}^0(x), \forall x \in \Omega. \quad (39)$$

In order to have $\Delta\hat{V}^i(x) < 0$ one needs $2c/(1-c)\overline{V}^0(x) < U(x, \hat{h}^i(x)) + \epsilon^i(x), \forall x$, which holds if

$$\frac{2c}{1-c}\overline{V}^0(x) < (1-c)U(x,0) \Leftrightarrow \frac{2c}{(1-c)^2} < \frac{U(x,0)}{\overline{V}^0(x)}. \quad (40)$$

Note that $\overline{V}^0(x) \le 2V^0(x), \forall x$, by definition of $\overline{V}^0(x)$ which is the value function of $h^{-1}(\cdot)$ with the utility of $U(x_k, u_k) + cU(x_k, 0)$, while, $V^0(x)$ is the value function of the same control policy with the utility of $U(x_k, u_k)$. Therefore, from $V^0(x) \le \gamma U(x,0)$ one has $\overline{V}^0(x) \le 2\gamma U(x,0)$. Hence, $1/(2\gamma) \le U(x,0)/\overline{V}^0(x), \forall x$. Therefore, if $2c/(1-c)^2 < 1/(2\gamma)$ or equivalently if

$$c^2 - (2+4\gamma)c + 1 > 0, \quad (41)$$

then inequality (40) holds. The root of the left hand side of the foregoing inequality are real and given by

$$c_1 = 1 + 2\gamma - \sqrt{4\gamma^2 + 4\gamma},\ c_2 = 1 + 2\gamma + \sqrt{4\gamma^2 + 4\gamma}. \quad (42)$$

Inequality (41) holds if $c < c_1$ or if $c > c_2$ by analysis of the sign of the quadratic equation on its left hand side. But, $c_2 > 1$, hence, any $c$ which satisfies $c > c_2$ will be unacceptable as such a $c$ does not belong to $[0,1)$. As for $c < c_1$ it is required to make sure $c_1 > 0$, otherwise no suitable $c$ will be resulted from this analysis, also. From $4\gamma^2 + 4\gamma + 1 > 4\gamma^2 + 4\gamma$, considering the non-negativeness of both sides, one has $\sqrt{4\gamma^2 + 4\gamma + 1} = 2\gamma + 1 > \sqrt{4\gamma^2 + 4\gamma}$, hence, $1 + 2\gamma - \sqrt{4\gamma^2 + 4\gamma} > 0$. Therefore, $c_1$ is indeed positive and a non-negative $c$ smaller than $c_1$ leads to the desired stability. The first part of the theorem is proved by noticing that when (30) holds, $\Delta\hat{V}^i(x)$ is strictly less than zero at any $x \ne 0$, considering the positive definiteness of $U(\cdot,0)$.

Inequality $\Delta\hat{V}^i(x) < 0$ leads to $\hat{\mathcal{B}}^i_{\bar{r}^i} \subset \Omega$ being an EROA for the system operated with $\hat{h}^i(\cdot)$. The reason is, any trajectory initiated within $\hat{\mathcal{B}}^i_{\bar{r}^i}$ will remain inside the set and hence, within $\Omega$ and therefore, converge to the origin. Finally, since $\hat{\mathcal{B}}^i_{\bar{r}^i}$ is contained in $\Omega$, it is bounded. Also, the set is closed, because, it is the *inverse image* of a closed set, namely $[0, \bar{r}^i]$ under a continuous function (by the continuity of function approximator,) [17]. Hence, $\hat{\mathcal{B}}^i_{\bar{r}^i}$ is compact. The origin is an *interior* point of the EROA, because $\hat{V}^i(0) = 0$, $\bar{r}^i > 0$, and $\hat{V}^i(\cdot)$ is continuous in $\Omega$. This completes the proof. $\square$

Theorem 2 proves that each single $\hat{h}^i(\cdot)$ if *constantly* applied on the system, will steer the states toward the origin. However, in online learning, the control policy will be subject to adaptation. In other words, if $h^i(\cdot)$ is applied at the current time, control policy $h^{i+1}(\cdot)$ might be the one which will be applied later. Therefore, another stability analysis is required for the system operated under *evolving/time-varying* policies. This is done for the general case of applying each policy $\hat{h}^i(\cdot)$ for $M_i \in \mathbb{N}$ steps before switching to the next version, i.e., $\hat{h}^{i+1}(\cdot)$ (and applying it for $M_{i+1} \in \mathbb{N}$ steps).

**Theorem 3.** *Let Assumption 2 hold and the system be operated using sequence of control policies $\{\hat{h}^i(\cdot)\}_{i=0}^\infty$, where each $\hat{h}^i(\cdot)$ is applied for $M_i \in \mathbb{N}$ time steps. Every trajectory contained in $\Omega$ will converge to the origin if $c$ is such that*

$$0 \le c < 1 + 4\gamma - \sqrt{16\gamma^2 + 8\gamma}, \quad (43)$$

*for some $\gamma \in \mathbb{R}_+$ that satisfies $V^0(x) \le \gamma U(x,0), \forall x \in \Omega$.*

*Proof*: Let the state trajectory generated through the scenario of applying each $\hat{h}^i(\cdot)$ at $M_i$ steps be denoted with $\hat{x}_k^+$. Eq. (12) and Lemma 1 lead to

$$\hat{V}^1(\hat{x}_0^+) = U(\hat{x}_0^+, \hat{h}^0(\hat{x}_0^+)) + \hat{V}^0(\hat{x}_1^+) + \epsilon^0(\hat{x}_0^+) \le \hat{V}^0(\hat{x}_0^+) + 2c\breve{V}^0(\hat{x}_0^+), \forall \hat{x}_0^+ \in \Omega. \quad (44)$$

Hence,
$$U(\hat{x}_0^+, \hat{h}^0(\hat{x}_0^+)) + \hat{V}^0(\hat{x}_1^+) + \epsilon^0(\hat{x}_0^+) - 2c\check{V}^0(\hat{x}_0^+) \leq \\ \hat{V}^0(\hat{x}_0^+), \forall \hat{x}_0^+ \in \Omega. \quad (45)$$

Evaluating (45) at $\hat{x}_1^+$ gives
$$U(\hat{x}_1^+, \hat{h}^0(\hat{x}_1^+)) + \hat{V}^0(\hat{x}_2^+) + \epsilon^0(\hat{x}_1^+) - 2c\check{V}^0(\hat{x}_1^+) \leq \\ \hat{V}^0(\hat{x}_1^+), \forall \hat{x}_1^+ \in \Omega. \quad (46)$$

Replacing $\hat{V}^0(\hat{x}_1^+)$ in (45) by the left hand side of (46) gives
$$U(\hat{x}_0^+, \hat{h}^0(\hat{x}_0^+)) + U(\hat{x}_1^+, \hat{h}^0(\hat{x}_1^+)) + \hat{V}^0(\hat{x}_2^+) + \\ \epsilon^0(\hat{x}_0^+) + \epsilon^0(\hat{x}_1^+) - 2c\check{V}^0(\hat{x}_0^+) - 2c\check{V}^0(\hat{x}_1^+) \leq \\ \hat{V}^0(\hat{x}_0^+), \forall \hat{x}_0^+ \in \Omega. \quad (47)$$

This process, i.e., using (45) 'in itself', may be repeated for the total of $M_0 - 1$ times to get
$$\sum_{k=0}^{M_0-1} \left( U(\hat{x}_k^+, \hat{h}^0(\hat{x}_k^+)) + \epsilon^0(\hat{x}_k^+) - 2c\check{V}^0(\hat{x}_k^+) \right) \\ + \hat{V}^0(\hat{x}_{M_0}^+) \leq \hat{V}^0(\hat{x}_0^+), \forall \hat{x}_0^+ \in \Omega. \quad (48)$$

Similarly, from Eq. (12) and Lemma 1, one has
$$\hat{V}^2(x) = U(x, \hat{h}^1(x)) + \hat{V}^1(f(x, \hat{h}^1(x))) + \epsilon^1(x) \leq \\ \hat{V}^1(x) + 2c\check{V}^1(x), \forall x \in \Omega. \quad (49)$$

Bringing $2c\check{V}^1(x)$ to the left side of the inequality and evaluating the result at $\hat{x}_{M_0}^+$ and $\hat{x}_{M_0+1}$, respectively, lead to
$$U(\hat{x}_{M_0}^+, \hat{h}^1(\hat{x}_{M_0}^+)) + \hat{V}^1(\hat{x}_{M_0+1}^+) + \epsilon^1(\hat{x}_{M_0}^+) \\ - 2c\check{V}^1(\hat{x}_{M_0}^+) \leq \hat{V}^1(\hat{x}_{M_0}^+), \forall \hat{x}_{M_0}^+ \in \Omega. \quad (50)$$

and
$$U(\hat{x}_{M_0+1}^+, \hat{h}^1(\hat{x}_{M_0+1}^+)) + \hat{V}^1(\hat{x}_{M_0+2}^+) + \epsilon^1(\hat{x}_{M_0+1}^+) \\ - 2c\check{V}^1(\hat{x}_{M_0+1}^+) \leq \hat{V}^1(\hat{x}_{M_0+1}^+), \forall \hat{x}_{M_0+1}^+ \in \Omega. \quad (51)$$

Replacing $\hat{V}^1(\hat{x}_{M_0+1}^+)$ in (50) using the left hand side of (51) leads to
$$U(\hat{x}_{M_0}^+, \hat{h}^1(\hat{x}_{M_0}^+)) + U(\hat{x}_{M_0+1}^+, \hat{h}^1(\hat{x}_{M_0+1}^+)) + \hat{V}^1(\hat{x}_{M_0+2}^+) \\ + \epsilon^1(\hat{x}_{M_0}^+) + \epsilon^1(\hat{x}_{M_0+1}^+) - 2c\check{V}^1(\hat{x}_{M_0}^+) - 2c\check{V}^1(\hat{x}_{M_0+1}^+) \\ \leq \hat{V}^1(\hat{x}_{M_0}^+), \forall \hat{x}_{M_0}^+ \in \Omega. \quad (52)$$

Repeating this process for $M_1 - 2$ more times results in
$$\sum_{k=0}^{M_1-1} \left( U(\hat{x}_{M_0+k}^+, \hat{h}^1(\hat{x}_{M_0+k}^+)) + \epsilon^1(\hat{x}_{M_0+k}^+) - 2c\check{V}^1(\hat{x}_{M_0+k}^+) \right) \\ + \hat{V}^1(\hat{x}_{M_0+M_1}^+) \leq \hat{V}^1(\hat{x}_{M_0}^+), \forall \hat{x}_{M_0}^+ \in \Omega. \quad (53)$$

Using inequality $\hat{V}^1(x) \leq \hat{V}^0(x) + 2c\check{V}^0(x)$, which is from Lemma 1, one can replace the $\hat{V}^1(\cdot)$ in the right hand side of (53) to get
$$\sum_{k=0}^{M_1-1} \left( U(\hat{x}_{M_0+k}^+, \hat{h}^1(\hat{x}_{M_0+k}^+)) + \epsilon^1(\hat{x}_{M_0+k}^+) - 2c\check{V}^1(\hat{x}_{M_0+k}^+) \right) \\ + \hat{V}^1(\hat{x}_{M_0+M_1}^+) - 2c\check{V}^0(\hat{x}_{M_0}^+) \leq \hat{V}^0(\hat{x}_{M_0}^+), \forall \hat{x}_{M_0}^+ \in \Omega. \quad (54)$$

Now, the left hand side of (54) can be used to replace $\hat{V}^0(x_{M_0}^+)$ in (48). This leads to
$$\sum_{k=0}^{M_0-1} \left( U(\hat{x}_k^+, \hat{h}^0(\hat{x}_k^+)) + \epsilon^0(\hat{x}_k^+) - 2c\check{V}^0(\hat{x}_k^+) \right) \\ + \sum_{k=0}^{M_1-1} \left( U(\hat{x}_{M_0+k}^+, \hat{h}^1(\hat{x}_{M_0+k}^+)) + \epsilon^1(\hat{x}_{M_0+k}^+) \\ - 2c\check{V}^1(\hat{x}_{M_0+k}^+) \right) + \hat{V}^1(\hat{x}_{M_0+M_1}^+) - 2c\check{V}^0(\hat{x}_{M_0}^+) \leq \\ \hat{V}^0(\hat{x}_0^+), \forall \hat{x}_0^+ \in \Omega. \quad (55)$$

The aforementioned inequality is resulted through incorporating control policies $\hat{h}^0(.)$ and $\hat{h}^1(.)$ applied for $M_0$ and $M_1$ steps, respectively, i.e., two generations of policies. This process can be repeated to include N generations of $\hat{h}^i(\cdot)$, leading to
$$\sum_{i=0}^{N-1} \sum_{k=0}^{M_i-1} \left( U(\hat{x}_{\sum_{j=0}^{i-1} M_j + k}^+, \hat{h}^i(\hat{x}_{\sum_{j=0}^{i-1} M_j + k}^+)) + \\ \epsilon^i(\hat{x}_{\sum_{j=0}^{i-1} M_j + k}^+) - 2c\check{V}^i(\hat{x}_{\sum_{j=0}^{i-1} M_j + k}^+) \right) \\ - \sum_{i=1}^{N-1} 2c\check{V}^i(\hat{x}_{\sum_{j=0}^{i-1} M_j}^+) + \hat{V}^{N-1}(\hat{x}_{\sum_{j=0}^{N-1} M_j}^+) \leq \\ \hat{V}^0(\hat{x}_0^+), \forall \hat{x}_0^+ \in \Omega. \quad (56)$$

On another hand, one has
$$2c\check{V}^0(\hat{x}_0^+) + \sum_{i=0}^{N-1} \sum_{k=1}^{M_i-1} 2c\check{V}^i(\hat{x}_{\sum_{j=0}^{i-1} M_j + k}^+), \quad (57) \\ \geq 0, \forall \hat{x}_0^+ \in \Omega,$$

as all the terms in the left hand side of the foregoing inequality are non-negative. Subtracting this left hand side from the left hand side of (56) leads to
$$\sum_{i=0}^{N-1} \sum_{k=0}^{M_i-1} \left( U(\hat{x}_{\sum_{j=0}^{i-1} M_j + k}^+, \hat{h}^i(\hat{x}_{\sum_{j=0}^{i-1} M_j + k}^+)) + \\ \epsilon^i(\hat{x}_{\sum_{j=0}^{i-1} M_j + k}^+) - 4c\check{V}^i(\hat{x}_{\sum_{j=0}^{i-1} M_j + k}^+) \right) \\ + \hat{V}^{N-1}(\hat{x}_{\sum_{j=0}^{N-1} M_j}^+) \leq \hat{V}^0(\hat{x}_0^+), \forall \hat{x}_0^+ \in \Omega. \quad (58)$$

From (37) and $\overline{V}^0(x) \leq 2V^0(x) \leq 2\gamma U(x, 0), \forall x \in \Omega$, one has
$$\check{V}^i(x) \leq \frac{2\gamma}{1-c} U(x, 0), \forall x \in \Omega, \forall i \in \mathbb{N}. \quad (59)$$

Also, $(1-c)U(x,0) \leq U(x,u) + \epsilon^i(x), \forall x, \forall u, \forall i$. Therefore, from (58) one has
$$(1 - c - \frac{8c\gamma}{1-c}) \sum_{i=0}^{N-1} \sum_{k=0}^{M_i-1} U(\hat{x}_{\sum_{j=0}^{i-1} M_j + k}^+, 0) \\ + \hat{V}^{N-1}(\hat{x}_{\sum_{j=0}^{N-1} M_j}^+) \leq \hat{V}^0(\hat{x}_0^+), \forall \hat{x}_0^+. \quad (60)$$

Assume
$$(1-c) - \frac{8c\gamma}{1-c} > 0 \Leftrightarrow c^2 - (2+8\gamma)c + 1 > 0. \quad (61)$$

Considering (60), the desired stability result can be obtained for $\hat{x}_k^+$s, providing (61) holds. That is, the sequence of partial sums of the left hand side of (60) is upper bounded and because of being non-decreasing it converges, as $N \to \infty$, [17], leading to the convergence of states to the origin, as long as the entire state trajectory is contained in $\Omega$. Finally, in order to enforce (61), one will need (43). The details are identical to the respective part in proof of Theorem 2, as replacing $\gamma$ in (41) with $2\gamma$, the left hand sides of inequalities (41) and (61) become identical. Moreover, the proof of positiveness of the right hand side of (43) follows from the same argument presented in that proof. $\square$

Finally, the last step of our analysis is presenting some results regarding the EROA for the system operated using evolving control policy during AVI.

**Theorem 4.** *Let conditions of Theorem 3 hold. Moreover, let $\hat{\mathcal{B}}_r^i := \{x \in \mathbb{R}^n : \hat{V}^i(x) \leq r\}$ and $\mathcal{B}_r^* := \{x \in \mathbb{R}^n : V^*(x) \leq r\}$ for any $r \in \mathbb{R}_+$. If $\mathcal{B}_r^* \subset \Omega$ for a given $r > 0$, then compact set $\hat{\mathcal{B}}_{(1-c)r}^0$ is an estimation of the region of attraction of the system under the evolving control policy of applying each $\hat{h}^i(\cdot)$ for $M_i \in \mathbb{N}$ time steps.*

*Proof*: Once conditions of Theorem 3 hold, (60) gives
$$\hat{V}^{N-1}(\hat{x}^+_{\sum_{j=0}^{N-1} M_j}) \leq \hat{V}^0(\hat{x}^+_0), \forall \hat{x}^+_0 \in \Omega, \forall N \in \mathbb{N} - \{0\}, \quad (62)$$
where $\hat{x}^+_k$ denotes the state trajectory generated by applying each $\hat{h}^i(\cdot)$ for $M_i \in \mathbb{N}$ steps. Therefore,
$$x^+_0 \in \hat{\mathcal{B}}_r^0 \Rightarrow \hat{x}^+_{\sum_{j=0}^{N-1} M_j} \in \hat{\mathcal{B}}_r^{N-1}, \forall N \in \mathbb{N} - \{0\}, \forall r. \quad (63)$$
From $\underline{V}^k(x) \leq \hat{V}^k(x), \forall x, \forall k$, cf. Theorem 1, one has
$$\hat{\mathcal{B}}_r^k \subset \underline{\mathcal{B}}_r^k, \forall k \in \mathbb{N}, \forall r \in \mathbb{R}_+, \quad (64)$$
where $\underline{\mathcal{B}}_r^i := \{x \in \mathbb{R}^n : \underline{V}^i(x) \leq r\}$ and $\underline{V}^k(x)$ is defined in Theorem 1. Moreover, $\{\underline{V}^k(x)\}_{k=0}^\infty$ is non-increasing and converges to $\underline{V}^*(x)$, [2], i.e., the optimal value function corresponding to cost function (18), because it is resulted from an exact VI. Therefore, $\underline{V}^*(x) \leq \underline{V}^{k+1}(x) \leq \underline{V}^k(x), \forall x$, which leads to
$$\underline{\mathcal{B}}_r^k \subset \underline{\mathcal{B}}_r^{k+1} \subset \underline{\mathcal{B}}_r^*, \forall k \in \mathbb{N}, \forall r \in \mathbb{R}_+. \quad (65)$$
From (64) and (65) one has
$$\hat{\mathcal{B}}_r^k \subset \underline{\mathcal{B}}_r^*, \forall k \in \mathbb{N}, \forall r \in \mathbb{R}_+. \quad (66)$$
On the other hand, it can be shown that $V^*(x) \leq (1-c)^{-1} \underline{V}^*(x), \forall x$, which if holds leads to
$$\underline{\mathcal{B}}_r^* \subset \mathcal{B}_{r/(1-c)}^*, \forall r \in \mathbb{R}_+. \quad (67)$$
Let the state trajectory generated from the exact VI of $\underline{V}^*(x)$ be denoted with $\underline{x}_k^*, \forall k$. Comparing cost function (18) with (2) one has
$$V^*(x_0) \leq \underline{V}^*(x_0) + c \sum_{k=0}^\infty U(\underline{x}_k^*, 0), \forall x_0. \quad (68)$$
The reason is, both sides of the foregoing inequality are infinite sums over $U(\cdot, \cdot)$s, except that the the left hand side is evaluated along the optimal trajectory with respect to (2) while the right hand side is evaluated along the optimal trajectory with respect to (18). Moreover, by composition $(1-c) \sum_{k=0}^\infty U(\underline{x}_k^*, 0) \leq \underline{V}^*(x_0)$ which once used in (68) leads to $V^*(x) \leq (1-c)^{-1} \underline{V}^*(x)$.

Finally, from (63), (66), and (67) one has
$$x^+_0 \in \hat{\mathcal{B}}_r^0 \Rightarrow \hat{x}^+_{\sum_{j=0}^{N-1} M_j} \in \mathcal{B}^*_{r/(1-c)}, \forall N \in \mathbb{N} - \{0\}, \forall r. \quad (69)$$

Based on the foregoing result, if $x_0^+ \in \hat{\mathcal{B}}_r^0$, then, the state vector *at the end of* each period of applying a single $\hat{h}^i(.)$, denoted with $\hat{x}^+_{\sum_{j=0}^{N-1} M_j}$ is inside $\mathcal{B}^*_{r/(1-c)}$. We, however, are interested in showing that the *entire* state trajectory stays within $\mathcal{B}^*_{r/(1-c)}$. This is not difficult to establish, as (69) holds for any $M_{N-1}$. Therefore, if placement of the state vector within $\mathcal{B}^*_{r/(1-c)}$ at a given $k$ is of interest where $k$ is such that it corresponds to a time instant in the middle of the period of applying $\hat{h}^{N-1}(.)$, i.e., if $\sum_{j=0}^{N-2} M_j < k < \sum_{j=0}^{N-1} M_j$, selecting $M_{N-1} = \sum_{j=0}^{N-2} M_j - k$ one has $\sum_{j=0}^{N-1} M_j = k$. Therefore, the given state will be within $\mathcal{B}^*_{r/(1-c)}$ per (69). The reason for being able to adjust $M_{N-1}$ is, relation (69) was obtained using (60). The latter was simply obtained using
$$U(x, \hat{h}^{N-1}(x)) + \hat{V}^{N-1}(f(x, \hat{h}^{N-1}(x))) \\ + \epsilon^{N-1}(x) \leq \hat{V}^{N-1}(x) - 2c\breve{V}^{N-1}(x), \quad (70)$$
(which corresponds to control policy $\hat{h}^{N-1}(\cdot)$) 'in itself' for $M_{N-1}$ time steps and combining the result with similar results for previous generations of control policies. One can always do this process of using (70) in itself for a smaller $M_{N-1}$. It should be noted that the length of period of applying previous control policies, i.e., $M_j, j \leq N-2$, *cannot* be altered, as it could lead to deviation of the resulting state trajectory from the one denoted with $\hat{x}_k^+$. But, the *stopping point* for using (70) in itself (see proof of Theorem 3) for obtaining the desired equation to serve the purpose in here, can be adjusted. Therefore, (69) leads to
$$x_0^+ \in \hat{\mathcal{B}}_r^0 \Rightarrow \hat{x}_k^+ \in \mathcal{B}^*_{r/(1-c)}, \forall k \in \mathbb{N} - \{0\}, \forall r \in \mathbb{R}_+. \quad (71)$$
Hence, any state trajectory initiated within $\hat{\mathcal{B}}^0_{(1-c)r}$ remains inside $\mathcal{B}_r^*$ and if $r$ is such that $\mathcal{B}_r^* \subset \Omega$, the trajectory remains within $\Omega$ and by Theorem 3 converges to the origin. $\square$

## V. NUMERICAL ANALYSIS

Part of the theoretical results presented in this study is numerically verified in this section. Towards this goal, Van der Pol's oscillator, with continuous-time dynamics of $\ddot{z} = (1 - z^2)\dot{z} - z + u$ is selected. The problem was taken into state space by defining $x = [X, Y]^T := [z, \dot{z}]^T$ and discretized with sampling time $\Delta t = 0.05s$. Moreover, cost function terms $Q(x) = 0.5x^T x$, $R = 0.05$, and $U(x_k, u_k) := Q(x_k) + u_k^T R u_k$ were selected in (2). For implementation of the AVI, the initial admissible policy was selected as (feedback linearization based) policy $h^{-1}(x) = -(1-X^2)Y - 4X - 5Y$.

The function approximator was selected in a polynomial form made of elements of $x$ up to the fourth order. The region of interest was selected as $\Omega := [-1, 1] \times [-1, 1] \subset \mathbb{R}^2$. 200 random $x$s were selected from $\Omega$ in each evaluation of Eq. (8), which due to the presence of approximation errors leads to (10), and method of least squares was utilized for finding the parameters (coefficients of the polynomial terms).

Given the control affine nature of the system, the minimizer in (11) can be found by setting the gradient of the term subject to minimization to zero, leading to
$$u = -\frac{1}{2} R^{-1} g^T \nabla \hat{V}^i(f(x, u)), \quad (72)$$
where $\nabla \hat{V}^i(x) := (\partial \hat{V}(x)/\partial x)^T$ and $g := \Delta t[0, 1]^T$. Given

the point that the unknown $u$ exists on both sides of Eq. (72), the following successive approximation may be used for finding the unknown, [7].

$$u^{j+1} = -\frac{1}{2}R^{-1}g^T \nabla \hat{V}^i\big(f(x, u^j)\big), \quad (73)$$

Convergence tolerance of $10^{-4}$ was selected for evaluation of convergence of iterations in (73).

At each iteration of the training, the accuracy of the function approximation was evaluated using two different sets of sample states, namely, *training samples* and *test samples*. The training samples are the 200 randomly selected samples for conducting least squares. But, the test samples are *new* samples selected by gridding the state space using squares with width of 0.05, leading to 1681 samples. The test samples were not used in the training, so that one can evaluate the *generalization* capability of the function approximator. The approximation error $\epsilon^i(\cdot)$ was then numerically found through (10), once all other parameters are known. Given $\epsilon^i(\cdot)$ at different training and test samples, constant $c$ used in $|\epsilon^i(x)| \leq cU(x,0), \forall i, \forall x \in \Omega$ was numerically found by evaluating $\max |\epsilon^i(\cdot)|/U(\cdot,0)$ at different samples. This process led to $c = c_1 = 0.0020$ based on the training samples and $c = c_2 = 0.0021$ based on the test samples. The smaller $c_1$ compared with $c_2$ was expected, given the point that $c_1$ is based on the training samples. Interestingly, $c_2$ is only slightly greater than $c_1$, which demonstrates the good generalization capability of the selected function approximator.

Similarly, constant $\gamma$, used in $V^0(x) \leq \gamma U(x,0), \forall x \in \Omega$ was numerically calculated to be $\gamma = \gamma_1 = 24.07$ and $\gamma = \gamma_2 = 24.13$, for training and test samples, respectively. It may be mentioned that $\gamma$ should be selected based on $V^0(\cdot)$, not based on $\hat{V}^0(\cdot)$. But, the former is not available, as all we have is its approximation, given by the latter. Denoting the constant calculated based on $\hat{V}^0(\cdot)$ with $\hat{\gamma}$, i.e., $\hat{V}^0(x) \leq \hat{\gamma}U(x,0), \forall x$, constant $\gamma$ can be calculated using $\hat{\gamma}$. The reason is, one has $V^0(x) \leq 1/(1-c)\hat{V}^0(x)$ per (17), which along with $\hat{V}^0(x) \leq \hat{\gamma}U(x,0)$ leads to $V^0(x) \leq (\hat{\gamma}/(1-c))U(x,0)$. Therefore, $\gamma = \hat{\gamma}/(1-c)$. Given these calculations, the upper bound for $c$, given by (30) was found to be $0.010$. Considering $c_1 = 0.0020 < 0.010$ and $c_2 = 0.0021 < 0.010$, the results established in Theorem 2 hold.

Selecting the iteration index of $i = 6$, calculation of EROA denoted with $\hat{\mathcal{B}}_{\bar{r}}^i$ in Theorem 2 is the next step. Numerically it was found that $\bar{r} = 14.39$ is the greatest $r$ using which $\hat{\mathcal{B}}_r^i \subset \Omega$. Given this value for $\bar{r}$, region $\hat{\mathcal{B}}_{\bar{r}}^i$ is plotted in Fig. 1. Also, different initial conditions were selected and the respective state trajectories under the control policy calculated based on the sixth iteration of the value function are plotted in the same figure. It can be observed that the state trajectories did not leave the EROA and converged to the origin, as expected. These results confirm the ones given by Theorem 2.

## VI. Conclusions

Stability of a system under value iteration initiated using an admissible guess was established without ignoring approximation errors in the iterations. Straight forward conditions for guaranteeing the stability of the system under a fixed as well as an evolving control policy were developed. Moreover, regions of attraction for the two cases were established.

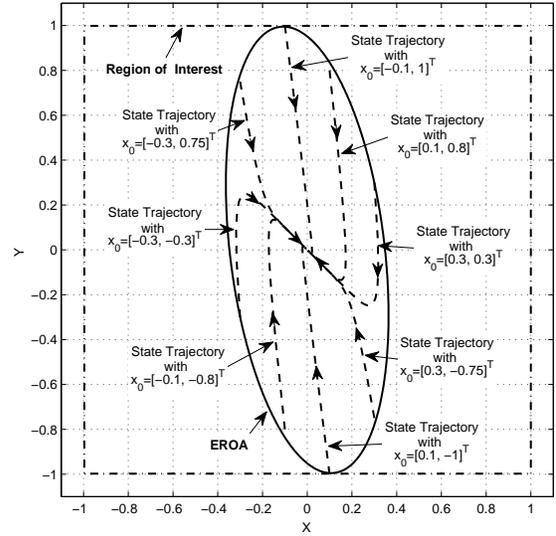

Fig. 1. State trajectories resulting from different initial conditions and the estimated region of attraction at $i = 6$.